\documentclass[12pt]{amsart}
\usepackage{graphicx}
\usepackage{amsmath}
\usepackage{amsthm}
\usepackage{amssymb}
\usepackage[textwidth=3.5cm]{todonotes}
\usepackage{tikz}
\usetikzlibrary{positioning}

\newtheorem{prob}{Problem}
\newtheorem{thm}{Theorem}
\newtheorem{prop}[thm]{Proposition}

\newtheorem{cor}[thm]{Corollary}

\newtheorem*{thm*}{Theorem}
\newtheorem*{lem*}{Lemma}

\theoremstyle{definition}
\newtheorem{defn}[thm]{Definition}
\newtheorem*{defn*}{Definition}

\newtheorem{ex}[thm]{Example}

\def\PA{\mathsf{PA}}

\def\lpa{{\mathcal L}_{\PA}}

\def\el{\mathcal L}
\def\elm{${\mathcal L}(\mc{M})$}
\def\m3{{\mathbf M_3}}
\def\mn{{\mathbf M_n}}
\def\n5{{\mathbf N_5}}
\def\b2{{\mathbf B_2}}
\def\bn{{\mathbf B_n}}
\def\lat{$\Lt(\mc{M})$}

\def\latp{$\Lt^+(\mc{M})$}
\def\Lat{$\Lt(\mc{N} / \mc{M})$}
\def\Lato{$\Lt_0(\mc{N} / \mc{M})$}
\def\Latp{$\Lt^+(\mc{N} / \mc{M})$}
\def\Latr{$\Ltr(\mc{N} / \mc{M})$}
\newcommand{\anglebracket}[1]{\langle #1 \rangle}

\newcommand{\mc}[1]{\mathcal{#1}}

\def\into{\longrightarrow}
\def\ct{countable}
\def\rs{recursively saturated}
\def\as{arithmetically saturated}
\def\st{structure}
\def\ar{arithmetic}
\def\vp{\varphi}
\def\om{\omega}

\DeclareMathOperator{\Eq}{Eq}
\DeclareMathOperator{\Aut}{Aut}
\DeclareMathOperator{\Cg}{Cg}
\DeclareMathOperator{\Lt}{Lt}
\DeclareMathOperator{\Ltr}{Ltr}
\DeclareMathOperator{\scl}{Scl}
\DeclareMathOperator{\Def}{Def}
\DeclareMathOperator{\Cod}{Cod}
\DeclareMathOperator{\SSy}{SSy}
\DeclareMathOperator{\Th}{Th}

\DeclareMathOperator{\Scl}{Scl}
\DeclareMathOperator{\tp}{tp}

\DeclareMathOperator{\eleme}{\prec_\text{\sf end}}
\DeclareMathOperator{\cof}{\prec_\text{\sf cof}}

\begin{document}

\title{The Lattice Problem for Models of $\PA$}
\author{Athar Abdul-Quader and Roman Kossak}
\begin{abstract}
The lattice problem for models of Peano Arithmetic ($\PA$) is to determine which lattices can be represented as lattices of elementary submodels of a model of $\PA$, or, in greater generality, for a given model $\mc{M}$, which lattices can be represented as  interstructure lattices of elementary submodels $\mc{K}$ of an elementary extension $\mc{N}$ such that $\mc{M}\preccurlyeq \mc{K}\preccurlyeq\mc{N}$. The problem has been studied for the last 60 years and the results and their proofs show an interesting interplay between the model theory of PA, Ramsey style combinatorics, lattice representation theory, and elementary number theory.  We present a survey of the most important results together with a detailed analysis of some special cases to explain and motivate a technique developed by James Schmerl for constructing elementary extensions with prescribed interstructure lattices. The last section is devoted to a discussion of lesser-known results about lattices of elementary submodels of countable recursively saturated models of PA.
\end{abstract}

\maketitle

\section{Introduction}
This paper concerns the study of substructure and interstructure lattices of models of $\PA$ (Peano Arithmetic). All models will be models of $\PA$. We use $\mc{M}$, $\mc{N}$, $\mc{K}$\dots for models, and $M$, $N$, $K$, \dots for their domains. We write $\mc{M}\preccurlyeq\mc{N}$ if $\mc{M}$ is an elementary submodel of $\mc{N}$ and $\mc{M}\prec\mc{N}$ if $\mc{M}$ is a proper elementary submodel of $\mc{N}$ (i.e., if $\mc{M}\preccurlyeq \mc{N}$ and $\mc{M}\not=\mc{N}$).

For a subset $X$ of the domain of a model $\mc{N}$, $\Scl^{\mc{N}}(X)$ denotes the Skolem closure of $X$ in $\mc{N}$. Thanks to the provability of the least number principle in $\PA$, the \emph{Skolem closure} of  $X$ in $\mc{N}$ is the same as its definable closure, i.e., the closure of $X$ under all parameter-free definable functions.
For $\mc{M}\prec \mc{N}$ and $a\in N\setminus M$, we will denote $\Scl^{\mc{N}}(M\cup\{a\})$ by $\mc{M}(a)$ and, for emphasis, we will denote $\Scl^{\mc{N}}(0)$ by $\mc{N}_{\rm min}$. $\mc{N}_{\rm min}$ is the prime model of the complete theory of $\mc{N}$, which we denote by $\Th(\mc{N})$. 

We will discuss two types of lattices, namely substructure lattices and interstructure lattices. 
Given $\mc{M}\preccurlyeq \mc{N}$, the
\emph{interstructure lattice}, denoted  $\Lt(\mc{N} / \mc{M})$, consists of all models $\mc{K}$ such that $\mc{M}\preccurlyeq\mc{K}\preccurlyeq\mc{N}$ ordered by inclusion. Clearly $\mc{M}$ and $\mc{N}$ are the minimum and the maximum elements of $\Lt(\mc{N} / \mc{M})$, respectively.

\Lat\ is a complete lattice. The meet of any set  of elementary submodels of $\mc{N}$ is the intersection of all the models in the set, and the join is the Skolem closure of their union.

The \emph{substructure lattice} of $\mc{N}$, denoted $\Lt(\mc{N})$, is  $\Lt(\mc{N} / \mc{N}_{\rm min})$, i.e., it is the lattice of all elementary submodels of $\mc{N}$.

 The systematic study of substructure and interstructure lattices of models of $\PA$ begins with Haim Gaifman's seminal paper \cite{gaifman76}.  Gaifman introduced a powerful technique of end-extensional and minimal types, which he applied  to construct models with some specific substructure and interstructure lattices. In particular, he proved that for every set $I$,  every model $\mc{M}$ has an elementary extension $\mc{N}$ such that $\Lt(\mc{N} / \mc{M})$ is isomorphic to the Boolean algebra of all subsets of $I$. 
Soon after \cite{gaifman76} appeared, Schmerl \cite{schmerl78} confirmed Gaifman's conjecture that for every finite distributive lattice $D$, every model $\mc{M}$
has an elementary end extension $\mc{N}$ such that $\Lt(\mc{N} / \mc{M})$ is isomorphic to $D$. Recall that a lattice $D$ is \emph{distributive} if the operations of join and meet in $D$ distribute over each other, i.e.,  $x\land (y\lor z)=(x\land y)\lor (x\land z)$ and $x\lor (y\land z)=(x\lor y)\land (x\lor z)$. By well-known representation theorems of Stone and Priestley, a lattice $D$ is distributive if for some set $X$, $D$ is isomorphic to a sublattice of the Boolean algebra of all subsets of $X$.

George Mills \cite{mills} extended Gaifman's technique to types with arbitrary sets of variables and completely characterized all distributive lattices that can be represented as $\Lt(\mc{N} / \mc{M})$.  Almost at the same time, Paris \cite{paris77} and Alex Wilkie \cite{wilkie} applied different techniques to give examples of finite nondistributive substructure and interstructure lattices. We will give precise statements of all these results in Section \ref{nondist}.

The work on  substructure and interstructure lattices of models of $\PA$ has turned out to be intimately connected with  problems in lattice theory, Ramsey style combinatorics, and even some elementary number theory. Chapter 4 of \cite{ks} gives a comprehensive account  of the main results obtained prior to 2005. Our goal is to give a brief survey of the area, to review the main lattice representation technique introduced by James Schmerl in \cite{schmerl_lat_86}, and to report on some more recent results. The paper ends with  a list of open problems. For now, let us just mention the most outstanding one.

\begin{prob}
Can every finite lattice be represented as $\Lt(\mc{N})$ for some model $\mc{N}$ of $\PA$? 
\end{prob}

The paper is organized as follows. Preliminaries in Section 2, along with introduction of the terminology and notation, contains a more detailed discussion of two particular examples that are included to illustrate  basic concepts and to motivate more general results that follow. 

The study of substructure and interstructure lattices of models of PA splits into two distinct areas: the distributive and nondistributive cases. This is reflected by splitting the brief overview of some major results into two  sections 3 and 4 in which each case is discussed separately. 

The next two sections are devoted  to an introduction of Schmerl's method of CPP-representations. In Section 5 we introduce basic definitions and illustrate them with a simple example. The purpose is not only to show how the method works, but also why it might be necessary. Section 6 proceeds in a similar fashion; a variant of the CPP-representations method is applied to lattices enriched by a rank function.

Section 7 is about substructure and interstructure lattices that are enriched further by adding the isomorphism relation between substructures. The idea leads to the concept of diverse models and diverse extensions that are illustrated by examples.

Section 8 is motivated by the interesting fact that the pentagon lattice can be realized by interstructure lattices of end extensions of countable models, but none of those extensions is conservative. This leads to the problem: what are possible sets of subsets of a model that are coded in elementary end extensions with prescribed interstructure lattices?

Section 9 is devoted to lattices of elementary submodels of \ct\ \rs\ models of $\PA$, and to one particular theorem that shows that for a large class of such models the isomorphism type of each model in the class is determined by its complete theory and its lattice of elementary substructures.

We end in Section 10 with some open questions.

\section{Preliminaries}
This section presents all basic notions and results  that will be needed for the discussion that follows. 
Some concepts are illustrated by simple examples and some preliminary results. A  full introduction to the model theory of $\PA$ is Richard Kaye's \cite{kaye} and all details directly related to the lattice problem can be found in \cite[Chapter 4]{ks}.

Let $\el$ be the language of $\PA$.
Given an expansion $\mc{L}^\prime$ of $\el$, $\PA^*$ is the theory consisting of the axioms of $\PA$ along with the induction schema for all formulas in $\mc{L}^\prime$. That is, $\PA^*$ is not a single theory, but many, one for each such expansion. With one notable exception \cite[Theorems 3 and 4]{schmerl_n5}, all results about models of $\PA$  in this article also hold for $\PA^*$ in any countable language. This also applies to the last section  in which we discuss \rs\ models, but there we need to add the assumption that the language is computable. 

A model $\mc{K}$ in \Lat\ is \emph{finitely generated} if $\mc{K}=\Scl^\mc{N}(M\cup A)$, for some finite $A\subseteq N$.
Each finite tuple in a model of $\PA$ is coded by a single element; hence, if $\mc{M}\prec \mc{N}$, for every finitely generated over $\mc{M}$ model $\mc{K}$ in \Lat\ there is an $a\in N$ such that $\mc{K}=\mc{M}(a)$.   

Let $L$ be a lattice. An element $x \in L$ is called \emph{compact} if whenever $x \leq \bigvee X$ for $X \subseteq L$, then $x \leq \bigvee X^\prime$ for some finite $X^\prime \subseteq X$. It is easy to verify that  the compact elements of \Lat\ are the finitely generated over $\mc{M}$ elementary submodels.

The set of compact elements of $\Lt(\mc{N} / \mc{M})$ forms a join-semilattice, which we denote as $\Lt_0(\mc{N} / \mc{M})$ and $\Lt_0(\mc{N})$ will denote $\Lt_0(\mc{N}/\mc{N}_{\rm{min}})$. For $a, b \in N$, $\mc{M}(a) \vee \mc{M}(b) = \mc{M}(\anglebracket{a, b})$, where $\anglebracket{x, y}$ is Cantor's pairing function. It is not obvious that the intersection of two finitely generated models may not be finitely generated. Various examples can be shown, one is given shortly below.

A complete lattice $L$ is \emph{algebraic} if each element of $L$ is the supremum of a set of compact elements. A lattice is $\kappa$-\emph{algebraic} if it is algebraic and each compact element has less than $\kappa$ many compact elements below it. 
The domain of any $\mc{K} \in \Lt(\mc{N} / \mc{M})$ is the union of the domains of all models $\mc{M}(a)$, for $a\in K$; thus $\mc{K}$ is the supremum  of the set of compact elements below it in the lattice. 

It follows from the remarks above that for any $\mc{M} \models \PA$, $\Lt(\mc{M})$ is $\aleph_1$-algebraic, and for any $\mc{M}$ and $\mc{N}$, if $\mc{M}\prec \mc{N}$ then   $\Lt(\mc{N} / \mc{M})$ is $|M|^+$-algebraic. Moreover,  if $\Lt_0(\mc{N}/ \mc{M}) \cong \Lt_0(\mc{N}_1 / \mc{M}_1)$, then $\Lt(\mc{N} / \mc{M}) \cong \Lt(\mc{N}_1 / \mc{M}_1)$. Thus, in order to realize a lattice as an interstructure lattice, we need only to ensure that we have control over the compact elements. Of course, if \Lat\ is finite, all elements of it are compact.

\subsection{Extensions}
A model $\mc{N}$ is a \emph{cofinal extension} of $\mc{M}$, if for every $b\in N$ there is an $a\in M$ such that $b<a$. $\mc{N}$ is an \emph{end extension} of $\mc{M}$ if for every $a\in M$ and $b\in N\setminus M$, $\mc{N} \models a < b$. We write $\mc{M} \cof \mc{N}$ if $\mc{N}$ is a cofinal elementary extension, and $\mc{M} \eleme \mc{N}$ if $\mc{N}$ is an elementary end extension.

If $\mc{M}\preccurlyeq\mc{N}$, then we say that $X\subseteq N$ is $\mc{M}$-\emph{definable} if it is defined in $\mc{N}$ by an $\mc{L}(\mc{M})$-formula. If $X\subseteq M$ is $\mc{M}$-definable, we just say that it is definable. 

 If $\mc{M}\preccurlyeq\mc{N}$, then $\mc{N}$ is a \emph{conservative} extension of $\mc{M}$, if for every definable $X\subseteq N$, $X\cap M$ is $\mc{M}$-definable. 

If $\mc{N}$ is a conservative extension of $\mc{M}$, then for each $c\in N$,  $A=\{x\in M: x<c\}$ is  $\mc{M}$-definable.  If $c<a$ for some $a\in M$,  then $A$, as a bounded definable subset of $\mc{M}$, must have a maximum element and it is easy to see that that maximum element is $c$; hence $c\in M$.  Therefore, all conservative extensions are end extensions. 

Let us note that $\Lt(\mc{N})$ may be uncountable even if $\mc{N}$ is finitely generated, i.e., it is finitely generated over $\mc{N}_{\rm min}$. To see this, let $\mc{M}$ be a countable model with an uncountable substructure lattice (as an example, take $\mc{M}$ to be a countable, recursively saturated model). By \cite[Theorem 2.1.12]{ks}, every countable model $\mc{M}$ has a \emph{superminimal} elementary end extension, i.e., an elementary end extension $\mc{N}$ such that $\mc{N}=\scl(b)$ for every $b\in N\setminus M$. 

For any completion $T$ of $\PA$, starting with a prime model of $T$, we can build an $\om_1$-chain of superminimal elementary  end extensions. The union of  this chain $\mc{N}$  has no proper elementary un\ct\ submodels and $\Lt_0(\mc{N})$ is isomorphic to $(\om_1,\leq)$. The existence of such models was proved independently  by Julia Knight \cite{Knight} and Jeff Paris \cite{paris72}.

\begin{ex}
We can use the technique of superminimal extension to provide an example showing that the intersection of two finitely generated submodels might not be finitely generated.

As mentioned above, every \ct\ model has a superminimal elementary end extension. Moreover, if  $\mc{M}$ is countable, nonstandard, and is generated by a bounded set of generators, then it has a superminimal cofinal extension (\cite[Exercise 2.5.2]{ks}). Suppose $\mc{M}$ is countable, nonstandard, and has a bounded, but not finite, set of generators. Such a model can be obtained as the union of an $\om$-chain of elementary cofinal extensions of a nonstandard finitely generated model $\mc{M}$, such that, for some fixed $a$ in $M$, each next model in the chain is a cofinal extension of the previous one generated over the previous model by a new element added below $a$. Let $\mc{N}_1$ be a superminimal elementary end extension of $\mc{M}$, $\mc{N}_2$ a superminimal cofinal extension of $\mc{M}$, and let $\mc{N} = \Scl(a, b)$, where $a \in N_1 \setminus M$ and $b \in N_2 \setminus M$ (by superminimality, any such $a$ and $b$ will work). Then, in $\Lt(N)$, $\mc{N}_1$ and $\mc{N}_2$ are compact (finitely generated), but their intersection is $\mc{M}$, not finitely generated.\end{ex}

\subsection{Ranked lattices}\label{ranked-def}
The following result is known as Gaifman's Splitting Theorem. 
\begin{thm}\label{gaifman} If $\mc{M}\preccurlyeq \mc{N}$, then there  is a unique ${\mc{K}}$ such that $\mc{M}\preccurlyeq_{\sf cof }{\mc{K}}\preccurlyeq_{\sf end}\mc{N}$. The domain of ${\mc{K}}$ is $\{x\in N: \exists y\in M\ \mc{N}\models (x\leq y)\}$. 
\end{thm}
Thanks to Gaifman's Splitting Theorem, for  each $\mc{K}$  in \Lat, we can define the \emph{rank} of $\mc{K}$,  $\rho(\mc{K})$, as the unique $\overline{\mc{K}}$ such that $\mc{K}\preccurlyeq_{\sf cof}\overline{\mc{K}}\preccurlyeq_{\sf end} \mc{N}$. The set of ranks of all models in \Lat, called the \emph{rankset}, is linearly ordered by inclusion. It is easy to see that for $\mc{K}_0$, $\mc{K}_1$ in \Lat, if $\mc{K}_0\prec \mc{K}_1$, then 
\begin{itemize}
\item $\mc{K}_0\preccurlyeq_{\sf cof}\mc{K}_1$ if and only if $\rho(\mc{K}_0)=\rho(\mc{K}_1)$, and \item $\mc{K}_0\preccurlyeq_{\sf end}\mc{K}_1$ if and only if $\mc{K}_0=\rho(K_0)\cap\mc{K}_1$.
\end{itemize}
Ranked lattices  were introduced by Schmerl in \cite{schmerl_lat_86}. Here is a modified definition from \cite{ks}: A \emph{ranked lattice} $(L,\rho)$ is a lattice $L$ equipped with a function $\rho:L\into L$ such that for all $x$ and $y$ in $L$ we have
\begin{enumerate}
    \item $x\leq \rho(x)$;
    \item $\rho(\rho(x))=\rho(x)$;
    \item $\rho(x)\leq \rho(y) \textup{ or } \rho(y)\leq \rho(x) $;
    \item $\rho(x\lor y)=\rho(x)\lor \rho(y)$.
\end{enumerate}

The \emph{rankset} of a ranked lattice $(L,\rho)$ is $\{\rho(x): x\in L\}$.

\Lat\ equipped with the rank function defined above is a ranked lattice, denoted by \Latr.

A lattice can have many expansions to a ranked lattice.  To be represented as rank functions in interstructure lattices of models of arithmetic those expansion have to satisfy certain additional conditions. 

Andreas Blass \cite{blass} showed that the intersection of two finitely generated cofinal submodels of a model $\mc{M}$ must be cofinal in $\mc{M}$.\footnote{Blass' theorem is about models of \emph{full \ar}, i.e., $\PA^*$ in the language with function and relation symbols for all functions and relations on $\om$, but its proof works for models of $\PA$ as well.}
Thus, if $\mc{M}\prec \mc{N}$, and $(L,\rho)$ is isomorphic to \Latr, then $(L,\rho)$ must satisfy \emph{the Blass Condition}: for all compact $x,y\in L$, if $\rho(x)=\rho(y)$, then $\rho(x)=\rho(x\land y)$.

Less perspicuous is \emph{the Gaifman Condition}: for all $x,y,z\in L$, if $x<y<x\lor z$, $z=\rho(z)$, and $x\land z=y\land z$, then $x=y$. See \cite[Proposition 4.2.12]{ks}. In the next section, we explain how these conditions are used to show that some finite lattices cannot  be represented by interstructure lattices given by end extensions.

\subsection{End-extensional and minimal types}
This subsection is a brief  summary of some of the results from \cite{gaifman76}. The results are interesting on their own.  They found many applications in the model theory of $\PA$ and the notion of definable type became standard in general model theory. Here we will just list some definitions  and results  that  play a major role in the solution to the lattice problem for distributive lattices that is the subject of the next section.

Let $T$ be a completion of $\PA$ and let $\mc{M}_T$ be the prime model of $T$. When we say that $p(x)$ is a type of $T$, we mean that $p(x)$ is in the language of $T$ and it is consistent with $T$. A type $p(x)$ of $T$ is \emph{unbounded} if $t<x$  is in $p(x)$ for each constant Skolem term $t$ of $T$. 

For a model $\mc{M}$, by $\mc{L}(\mc{M})$ we will denote  $\mathcal{L}$ with added constant symbols for all elements of $M$. An \elm-{\it type} is a type in \elm\ that is finitely realizable in $\mc{M}$. Thanks to the availability of a definable pairing function, for most of the results we will discuss, it is enough to consider 1-types.  

An \elm-type  is \emph{unbounded} if $(a<x)\in p(x)$ for all $a\in M$. So, for a completion $T$ of $\PA$, an $\el$-type $p(x)$ is unbounded if it is an unbounded ${\mathcal L}(\mc{M}_T)$-type.

For every model $\mc{M}\models T$, every unbounded type of $T$ extends to an unbounded \elm-type \cite[Proposition 2.7]{gaifman76}. 

\begin{defn}[\cite{gaifman76}] For a model $\mc{M}$, an \elm-type $p(x)$  is \emph{definable} if for every $\mathcal L$-formula $\vp(x,y)$   there exists an $\mathcal L$-formula $\sigma_\vp(y)$ such that for all $a\in M$,
\[\vp(x,a)\in p(x) \textup{ iff } \mc{M}\models\sigma_\vp(a). \]
A type of a completion $T$ is \emph{definable} if it is   a  definable ${\mathcal L}(\mc{M}_T)$-type.
\end{defn}
It is easy to see that $\mc{N}$ is a conservative extension of $\mc{M}$ if and only if  for every $a\in N\setminus M$, $\tp(a/M)$ is definable. 
If $p(x)$ is a complete \elm-type  then by $\mc{M}(p)$ we denote the unique up to isomorphism Skolem closure of $M\cup\{b\}$ in an elementary extension of $\mc{M}$ in which $b$ realizes $p(x)$.

\begin{defn}[\cite{gaifman76}]
   Let $T$ be a completion of $\PA$, and let $p(x)$ be a type of $T$. Then, 
   \begin{enumerate}
     \item $p(x)$ is \emph{end-extensional} if for every $\mc{M}\models T$ and every unbounded complete \elm-type $q(x)$,  if $p(x)\subseteq q(x)$, then $\mc{M}(q)$ is an end extension of $\mc{M}$.  
     \item $p(x)$ \emph{minimal} if for every  $q(x)$ as above, $\mc{M}(q)$ is a \emph{minimal extension} of $\mc{M}$, i.e., $\Lt(\mc{M}(q)/\mc{M})$ has exactly two elements: bottom $\mc{M}$ and top $\mc{M}(q)$.
   \end{enumerate}
\end{defn} 
Gaifman proved that every minimal type is end-extensional, that every end-extensional type is definable, and  that minimal and end-extensional types exist in abundance: for every completion $T$ of $\PA$ there are continuum many independent minimal types of $T$.\footnote{A type $p(x)$ depends on a type $q(x)$ if for some Skolem term $t(x)$, for all  formulas $\varphi(x)$ in $p(x)$, $\varphi(t(x))$ is in $q(x)$. Two types are independent if neither depends on the other.} If  $p$ and $q$ are such types, then for all $\mc{M}$, $\mc{M}(p)$ and $\mc{M}(q)$ are not isomorphic. 
From the existence of minimal types for every completion of $\PA$ it follows that every model of $\PA$ has a minimal (conservative) elementary end extension.

\begin{ex}\label{ex5} Let $p(x)$ be a minimal type of $T$. For $\mc{M}\models T$, let $\mc{N}=\mc{M}(a)(b)(c)=\mc{M}(\langle a,b,c\rangle)$, where $a, b, c$ all realize $p(x)$. Let $\mc{K}$ be the supremum of $\mc{M}(a)$ in $\mc{M}(\langle a,c\rangle)$.

\[\overbrace{-----)==a==)}^{\mc{K}}\underbrace{--b--)}_{\mc{M}(b)\setminus M}\overbrace{==c==)}^{\mc{M}(\langle a,c\rangle)\setminus { K}} \]
Then, it follows from the facts about  minimal types proved  in \cite{gaifman76} that \Lat$\cong ({\mathcal P}(\{a,b,c\}),\subseteq)$,  $\tp(\langle a,b\rangle)=\tp(\langle a,c\rangle)$, and both types are definable. The example shows that there are definable types which are not end-extensional. In particular, $\tp(\langle a,b\rangle)$ is such a type, because $\mc{M}(b)(\langle a,c\rangle)$ is not an end extension of $\mc{M}(b)$. 
\end{ex}

The  technique of minimal and end-extensional types uses infinitary combinatorics of unbounded definable sets in models of $\PA$ and does not apply to  cofinal extensions. Cofinal extensions are obtained by realizing bounded types for which combinatorial arguments about unbounded definable sets are replaced by their analogs involving bounded definable sets satisfying  suitable notions of largeness. This makes a difference. For example,  while we know that every model $\mc{M}$ has a minimal elementary end extension,  we can only prove that  that every \emph{countable} nonstandard model has a minimal cofinal extension. This was proved by Blass (see \cite[Corollary 2.1.6]{ks}). It is a long-standing open question whether every model of $\PA$ has a minimal cofinal extension. 

The example of minimal extensions shows that a lattice---in this example the two element lattice $L=\{{\mathbf 0_L},{\mathbf 1_L}\}$---can be realized as \Lat\ in two much different ways.  For a countable nonstandard $\mc{M}$, $\mc{N}$ can be either an elementary end extension or elementary cofinal extension (by Gaifman's splitting theorem, it cannot be a mixed extension). 

\section{Distributive lattices}\label{hist1}
Two of the main results about the lattice problem from Gaifman's \cite{gaifman76} are:
\begin{enumerate}
\item For any set $I$, every model $\mc{M}$ has a (conservative) elementary end extension $\mc{N}$ such that $\Lt(\mc{N}/\mc{M})$ is isomorphic to the Boolean algebra of all subsets of $I$ \cite[Theorem 4.10]{gaifman76}.
\item For every countable linearly ordered set $(J,<)$, every model $\mc{M}$ has a (conservative) elementary end extension $\mc{N}$ such that \Lato\ is isomorphic to $(J,<)$ \cite[Corollary 5.3]{gaifman76}.
\end{enumerate}
The proof of (1) uses an iteration of elementary end extensions generated by a minimal type, as we showed in Example \ref{ex5} above. The result in (2) is a corollary of a theorem about the existence of a particular end-extensional type that is used to get an extension of $\mc{N}$ in one step. The proof of that theorem is difficult and it takes 15 pages. Gaifman stresses that in both results one gets a uniform operator that produces the required extension $\mc{N}$ in a uniform way. Notice that in both cases the interstructure lattice is distributive.

Soon after \cite{gaifman76} appeared, Mills managed to generalize Gaifman's results by extending the notions of end-extensional and minimal types to types with infinitely many variables. His main result  is the following theorem.

\begin{thm}[\cite{mills}]\label{mills}
Let $D$ be a distributive lattice. Then the following are equivalent.
\begin{enumerate}
    \item There exists a model $\mc{M}$ such that $\Lt(\mc{M})$ is isomorphic to $D$.
    \item Every model  $\mc{M}$ has an elementary extension $\mc{N}$ such that  $\Lt(\mc{N}/\mc{M})$ is isomorphic to $D$.
    \item $D$ is $\aleph_1$-algebraic.
\end{enumerate}
\end{thm}
The equivalence of (1) and (3) was proved independently by Paris \cite{paris72}. The proof of Mills' theorem is difficult and it involves many technical details. 

In \cite{mills}, Mills follows the statement of Theorem \ref{mills} with an instructive example. Let $L$ be the ordered unit interval $([0,1],\leq)$. $L$ is a complete distributive lattice, but it is not $\aleph_1$-algebraic. Every element of $L$ is compact; hence for every nonzero element, there are continuum many compact elements below it. Let $L'$ be the closely related lattice of initial segments in the set of rational numbers in [0,1] ordered by inclusion. In $L'$, the compact elements are closed intervals $[0,p]$. Because each $[0, p]$ contains only countably many compact segments, by Mills' theorem, $L'$ is isomorphic to a substructure lattice.

All  results about end-extensional types and their applications to the lattice problem in \cite{gaifman76} hold for $\PA^*$ in a \ct\ language. Gaifman asked if they hold for uncountable languages as well. In response, Mills gave a construction of a model of $\PA^*$ in a language with $\aleph_1$ function symbols that has no elementary end extension \cite{mills_noend}.

In terminology of Mills \cite{mills}, a definable type $p(x)$ of $T$  \emph{produces} a lattice $L$ if for every model $\mc{M}$ of $T$, $\Lt(\mc{M}(p)/\mc{M})$ is isomorphic to $L$. For example, minimal types produce the two element lattice. It is a special feature of Gaifman's technique that his special types produce prescribed lattices \emph{for all} models of $T$, not just the countable ones. This is reminiscent of the MacDowell-Specker theorem, which says that \emph{every} model of $\PA$ has an elementary end extension. The proof of this fact for countable models follows by a relatively straightforward omitting types argument. The types  that  produce $\mc{M}(p)$ given $\mc{M}$ and $p(x)$ do it for all models $\mc{M}$ regardless of the cardinality of their domains. 

Gaifman conjectured that for every completion $T$ of $\PA$ and every finite distributive lattice that has a unique atom, there is an end-extensional type $p(x)$ of $T$ that produces $D$. The conjecture was confirmed by Schmerl \cite{schmerl78} and by a different construction by Mills \cite{mills}. Mills also proved the following variant of Theorem \ref{mills} that fully characterizes all distributive lattices that can be produced by end-extensional types. 

\begin{thm}
Let $T$ be a completion of $\PA^*$ in a \ct\ language and let $D$ be a distributive lattice. Then the following are equivalent:
\begin{enumerate}
    \item There is a definable (end-extensional) type $p(x)$ of $T$ which produces $D$.
    \item $D$ is $\aleph_1$-algebraic, $\bigvee D$ is compact (and for any nonzero $\alpha, \beta\in D$,  $\alpha\land\beta$ is nonzero). 
\end{enumerate}
\end{thm}

In conclusion, for every model $\mc{M}$ and every distributive lattice $D$ that can be realized as an interstructure lattice, $D$ can be  realized as $\Lt(\mc{N}/\mc{M})$, where $\mc{N}$ is an elementary conservative end extension of $\mc{M}$. This completely solves the general representation problem for the distributive lattices,  but one can still ask more specific questions about the extension $\mc{N}$ for which $\Lt(\mc{N}/\mc{M})$ is isomorphic to a given lattice $D$. In the distributive case this adds more depth to the subject. In the nondistributive case it turns out to be necessary, because there are  nondistributive lattices which cannot  be realized as $\Lt(\mc{N}/\mc{M})$, where $\mc{N}$ is an end, or even mixed,  extension of $\mc{M}$.

\section{Nondistributive lattices}\label{nondist}

By $\bn$ we denote the Boolean algebra of all subsets of an $n$-element set. For every $n$, $\bn$ is distributive  and so are all of its sublattices. The following  are examples of nondistsibutive lattices: for all $n>2$,  $\mn$ which is the lattice with  $n+2$ elements that has a top, a bottom, and $n$ incomparable elements in between, the pentagon lattice $\n5$, and the hexagon lattice $\mathbf H$. Dedekind proved that a lattice is distributive if and only if it has no sublattice that is isomorphic to either $\m3$ or $\n5$ (for a proof, see \cite[Theorem 102]{gratzer}).

\begin{figure}[ht]
\centering
\begin{tikzpicture}[x=.5cm, y=.5cm]
  \tikzstyle{every node}=[draw, circle]
  \node (one) at (-4, 2) {$1$};
  \node (a) at (-6, 0) {$a$};
  \node (b) at (-2, 0) {$b$};
  \node (zero) at (-4, -2) {$0$};
  \draw (zero) -- (a) -- (one) -- (b) -- (zero);

  \node (o) at (4, 2) {$1$};
  \node (aa) at (2, 0) {$a$};
  \node (bb) at (4, 0) {$b$};
  \node (cc) at (6, 0) {$c$};
  \node (z) at (4, -2) {$0$};
  \draw (z) -- (aa) -- (o) -- (bb) -- (z) -- (cc) -- (o);
\end{tikzpicture}
\caption{The lattices $\mathbf{B}_2$ and $\mathbf{M}_3$.}
\label{fig1}
\end{figure}
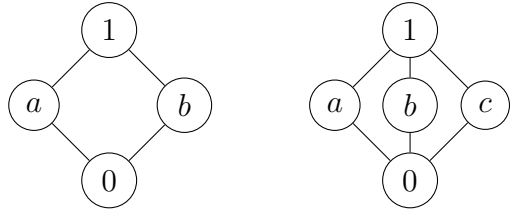

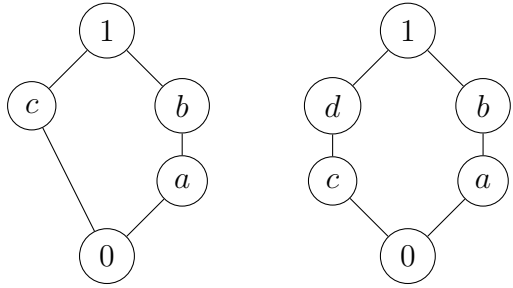
\begin{figure}[ht]
\centering
\begin{tikzpicture}[x=.5cm, y=.5cm]
  \tikzstyle{every node}=[draw, circle]
  \node (one) at (-4, 2) {$1$};
  \node (b) at (-2, 0) {$b$};
  \node (a) at (-2, -2) {$a$};
  \node (c) at (-6, 0) {$c$};
  \node (zero) at (-4, -4) {$0$};
  \draw (zero) -- (a) -- (b) -- (one) -- (c) -- (zero);

  \node (o) at (4, 2) {$1$};
  \node (bb) at (6, 0) {$b$};
  \node (aa) at (6, -2) {$a$};
  \node (dd) at (2, 0) {$d$};
  \node (cc) at (2, -2) {$c$};
  \node (z) at (4, -4) {$0$};
  \draw (z) -- (aa) -- (bb) -- (o) -- (dd) -- (cc) -- (z);
\end{tikzpicture}
\caption{The lattices $\mathbf{N}_5$ and $\mathbf{H}$.}
\label{fig2}
\end{figure}

All lattices $\bn$, as well as $\m3$,  $\n5$ and $\mathbf H$   can be represented as interstructure lattices, but the general results about how this is achieved are strikingly different. Here they are:
\begin{enumerate}
    \item (Gaifman \cite{gaifman76}) For each $n$, every model $\mc{M}$  has an elementary end extension $\mc{N}$ such that \Lat$\cong \bn$;
    \item (Wilkie \cite{wilkie}) Every \emph{countable} model $\mc{M}$ has an elementary end extension $\mc{N}$ such that \Lat$\cong {\n5}$.
     \item (Schmerl \cite{schmerl_lat_86}) Let  $L$ be either $\m3$ or $\mathbf H$. Then  every \emph{countable nonstandard} model  $\mc{M}$ has an elementary cofinal extension $\mc{N}$ such that \Lat$\cong L$.
\end{enumerate}
Independently, Gaifman \cite{gaifman76} and Paris \cite{paris72} showed that if $\mc{M}\prec_{\sf end} \mc{N}$ then \Lat\ is not isomorphic to $\m3$. Here is a short proof. Let $L=\m3$, and let $\rho$ be a rank on $L$ that satisfies the Blass Condition. We will show that  $\rho({\mathbf 0}_L)=\rho({\mathbf 1}_L)$. Suppose $a, b, c$ are the three incomparable elements of $L$. Because the rankset of a ranked lattice is linearly ordered, only one of these elements can be in the rankset of $(L,\rho)$. Suppose it is $c$. Then, by  condition (1) in the definition of rank, $\rho(a)=\rho(b)={\mathbf 1}_L$. By the Blass Condition, $\rho(a)=\rho(a\land b)=\rho({\mathbf 0}_L)$; hence $\rho({\mathbf 0}_L)=\rho({\mathbf 1}_L)$.
Thus, if $\m3$ is represented as \Lat, then $\mc{M}\prec_{\sf cof} \mc{N}$. 

Wilkie \cite{wilkie} proved that if $\mc{M}\prec_{\sf end} \mc{N}$ then \Lat\ is not isomorphic to the hexagon lattice $\mathbf H$. A short argument, due to Schmerl, using both the Blass and Gaifman conditions is given in \cite[Proposition 4.2.13]{ks}.

It is shown in \cite[Theorem 4.6.5]{ks} that if \Lat\ is isomorphic to $\n5$, then $\mc{N}$ is not a conservative extension of $\mc{M}$. Because there are uncountable models all of whose elementary end extensions are conservative (the \emph{rather classless models}), this shows that Wilkie's theorem about $\n5$ cannot be generalized to the uncountable case.

In the positive direction, much of what is known today about lattices that can be represented as \Lat\ rests on applications of  a powerful technique of representations of lattices introduced by Schmerl in \cite{schmerl_lat_86}. Section \ref{reps} is devoted to such representations. Here let us just note that for all lattices ${\mathbf M}_{q+1}$, where $q=p^k$ for a prime $p$ and $1\leq k$,   every nonstandard countable model $\mc{M}$ has a cofinal extension $\mc{N}$ such that $L$ is isomorphic to \Lat. 
Using this and other results about finite lattice representations, it can be shown that the smallest number $n$ for which $\mn$ is not known to have an interstructure lattice representation is 16 \cite[page 134]{ks}. 

The results above show that certain finite lattices can only be represented as \Lat\ only when $\mc{N}$ is a cofinal extension of $\mc{M}$. The following theorem shows that there are no finite lattices that can only be represented by interstructure lattices of end extensions.

\begin{thm}[\cite{schmerl_lat_86}]\label{cof} If a finite lattice $L$ can be represented as \Lat, for some $\mc{M}$ and $\mc{N}$, then every countable nonstandard model $\mc{M}'$ that is elementarily equivalent to $\mc{M}$ has a cofinal extension $\mc{N}'$ such that $\Lt(\mc{N}'/\mc{M}')$ is isomorphic to $L$.
\end{thm}

\subsection{The mysterious $\n5$} 
Let us take a closer look at representations of $\n5$ as interstructure lattices.

From Wilkie's and Schmerl's theorems it follows that every nonstandard \ct\ model $\mc{M}$ has a cofinal extension extension $\mc{N}$ such that \Lat\ is isomorphic to $\n5$. 

Suppose now that \Lat\ is a representation of $\n5$ and $\mc{M}$ is not cofinal in $\mc{N}$.
Let $F:\n5\into$ \Lat\ be an isomorphism. Using the labeling in Figure 2, let $F(r)=\mc{M}_r$, and let $\rho$ be a rank function of $\n5$.
Applying the Gaifman condition with $x=a$ and $y=b$ we get that $\rho(c)=1$, hence $\mc{M}_c\cof\mc{N}$.

Because $\rho(0)<1$ and $c\land b=0$, by the Blass condition $\rho(b)=b$. Hence $\mc{M}_b\eleme \mc{N}$. Finally, by \cite[Theorem 4.6.1]{ks}, $\rho(0)\not=b$, and we get that $\rho(0)$ must be either $0$ or $a$. In the first case we get that either $\mc{M}\eleme\mc{M}_a\eleme \mc{M}_b$ or $\mc{M}\eleme\mc{M}_a\cof\mc{M}_b$. As observed by Schmerl in \cite{schmerl_n5}, Wilkie's proof shows that both these scenarios can be realized. If $\rho(0)=a$, $\mc{N}$ is a mixed extension of $\mc{M}$. This case is surprisingly ruled out by the main result of \cite{schmerl_n5}. Thus $\mc{N}$ is not a mixed extension of $\mc{M}$.\footnote{Section 4.6 of \cite{ks} is devoted to representations of $\n5$. The claim there that it can be shown that $\n5$ can be realized as an interstructure lattice of a mixed extension turned out to be erroneous.} 

For $\mc{M}\prec\mc{N}$, $\Cod(\mc{N} / \mc{M})$ is the set of all intersections $M\cap X$, where $X$ ranges over all definable subsets of ${N}$.  In the short proof of \cite[Theorem 4.6.5]{ks}, a contradiction is derived from the assumptions that the extension $\mc{M}\eleme\mc{N}$ is conservative and that \Lat\ is isomorphic to $\n5$, but we do  not get any information about the undefinable subsets of $M$ that end up in $\Cod(\mc{N}/\mc{M}$).

In a recent paper \cite{schmerl_n5}, Schmerl adds more mystery the story  of representations of $\n5$. As we mentioned earlier, model-theoretic techniques  developed for models of $\PA$ often apply to models of $\PA^*$ in \ct\ languages. It turns out to be not so in the case of $\n5$. Schmerl proves that every \ct\ \rs\ model $\mc{M}$ of $\PA$ can be expanded to a model $\mc{M}^*$ of $\PA^*$ by adding countably many new sets so that $\mc{M}^*$ has a mixed elementary extension $\mc{N}^*$ such that $\Lt(\mc{N}^*/\mc{M}^*)$ is isomorphic to $\n5$. 

\subsection{Infinite nondistributive lattices}
In \cite{schmerl2010}, after brief remarks about the lattice problem for finite lattices, Schmerl writes ``Even less is known about infinite lattices; only in the case of distributive lattices had there been any significant results." Then he proceeds with an extended introduction to the main results in the paper that extend Theorem \ref{mills} and its earlier version dues to Paris \cite{paris72}, which---as corollaries---provide examples of infinite nondistributive interstructure lattices. In  particular the main theorem quoted below shows  that there is a model $\mc{M}$ for which $\Lt(\mc{M})$ is infinite and finitely generated. 

\begin{thm}[\cite{schmerl2010}] Let $L$ be an algebraic  bounded lattice with at most countably many compact elements. Then every countable nonstandard model $\mc{M}$ has a cofinal extension $\mc{N}$ such that \Lat$\cong L$.
\end{thm}

The theorems are formulated in lattice theoretic terms that take a good few pages to introduce.
We will only quote the following from \cite{day}. It  can be taken as a definition:  A finite lattice $L$ is \emph{bounded} if and only if it is in the smallest class containing the one element lattice and closed under the doubling of intervals. For example, $\n5$ is bounded, but $\m3$ is not.  A lattice is \emph{bounded} if each of its finitely generated sublattices is bounded. See \cite{schmerl2010} for definitions of all terms above.
 
The proofs heavily depend on major results in lattice representation theory that we will introduce in the next section.
 
\section{Representations}\label{reps}
Over the last 40 years, in several papers, Schmerl developed a special technique of constructing elementary extensions with prescribed interstructure lattices. It is based on particular representations of lattices as \emph{lattices of equivalence relations}. In this section, we give basic definitions and motivate them with examples.

\begin{defn} Let $A$ be any set and $L$ a finite lattice.
\begin{enumerate}
\item The set $\Eq(A)$ is the set of all equivalence relations on $A$. This set forms a lattice under inclusion, with $\mathbf{0}_A$ being the discrete relation $\{(a,a): a\in A\}$,   and $\mathbf{1}_A$ the trivial relation $A\times A$.
\item Let $L$ be a finite lattice and $A$ a set. Then $\alpha : L \into \Eq(A)$ is a \emph{pseudo-representation} of $L$ if:
\begin{itemize}
    \item $\alpha(0_L) = \mathbf{1}_A$, ($\alpha(0_L)$ is trivial)
    \item $\alpha(1_L) = \mathbf{0}_A$, ($\alpha(1_L)$ is discrete) and
    \item $\alpha(x \vee y) = \alpha(x) \wedge \alpha(y)$.
\end{itemize}
\item $\alpha$ is a \emph{representation} if it is a pseudo-representation and is one-to-one.
\end{enumerate}
\end{defn}

Oftentimes, realizing a particular finite lattice as an interstructure lattice requires choosing an appropriate representation of the lattice and proving some combinatorial lemmas about this representation. Before describing important properties of representations that are involved in these kinds of constructions, let us examine a motivating example.

\begin{ex}\label{b2}Let $\mc{M} \models \PA$ and $\mc{M} \prec \mc{N}$ such that $\Lt(\mc{N} / \mc{M}) \cong \mathbf{B}_2$. Then there are $a, b \in N$ such that $\mc{N} = \mc{M}(\anglebracket{a, b})$ and $\mc{M}(a)\not=\mc{M}(b)$ (see Figure 1). Let $\pi_1$ and $\pi_2$ be the projection functions, so that $$\PA \vdash \forall x \forall y (\pi_1(\anglebracket{x, y}) = x \wedge \pi_2(\anglebracket{x, y}) = y).$$ Then for any $c \in N$, there is an $\mc{M}$-definable $f : N \into N$ such that $\mc{N} \models f(\anglebracket{a, b}) = c$. Let us fix such  $c$ and $f$ and define a subset of $M$ as follows.   There are four possibilities:
\begin{itemize}
    \item $c \in M$. In this case, let $X = \{ \anglebracket{x, y} : \mc{M} \models f(\anglebracket{x, y}) = c \}$. Clearly, $\mc{M} \models \forall n, m \in X (f(n) = f(m))$.
    \item $c \in \mc{M}(a) \setminus \mc{M}$. In this case, there are $\mc{M}$-definable $g_1, g_2$ such that $$\mc{N} \models g_1(a) = c \wedge g_2(c) = a.$$ Let $X = \{ \anglebracket{x, y} : \mc{M} \models f(\anglebracket{x, y}) = g_1(x) \wedge g_2(f(\anglebracket{x, y})) = x \}$. Notice that $\mc{M} \models \forall n, m \in X (f(n) = f(m) \Longleftrightarrow \pi_1(n) = \pi_1(m))$.
    \item $c \in \mc{M}(b) \setminus \mc{M}$. Similarly, there are $\mc{M}$-definable $g_1, g_2$ such that $$\mc{N} \models g_1(b) = c \wedge g_2(c) = b.$$ Let $X = \{ \anglebracket{x, y} : \mc{M} \models f(\anglebracket{x, y}) = g_1(y) \wedge g_2(f(\anglebracket{x, y})) = y \}$. Again, one observes that $\mc{M} \models \forall n, m \in X (f(n) = f(m) \Longleftrightarrow \pi_2(n) = \pi_2(m))$.
    \item $c \in \mc{N} \setminus (\mc{M}(a) \cup \mc{M}(b))$. In this case, there is $\mc{M}$-definable $g$ such that $$\mc{N} \models g(c) = \anglebracket{a, b}.$$ Let $X = \{ \anglebracket{x, y} : \mc{M} \models g(f(\anglebracket{x, y})) = \anglebracket{x, y} \}$. Here we observe that $\mc{M} \models \forall n, m \in X (f(n) = f(m) \Longleftrightarrow n = m)$; that is, $f$ is one to one on $X$.
\end{itemize} In each of these cases, $X$ is an infinite, $\mc{M}$-definable, set and $\mc{N} \models \anglebracket{a, b} \in X$. Let us refer to an infinite, $\mc{M}$-definable $X$ such that $\anglebracket{a, b} \in X^{\mc{N}}$ as \emph{large}. 

Given any set $X$, we define the pseudo-representation $\alpha_X : \mathbf{B}_2 \into \Eq(X)$:
\begin{itemize}
    \item $\alpha_X(0)$ is trivial,
    \item $(n, m) \in \alpha_X(a)$ iff $\pi_1(n) = \pi_1(m)$,
    \item $(n, m) \in \alpha_X(b)$ iff $\pi_2(n) = \pi_2(m)$, and
    \item $\alpha_X(1)$ is discrete.
\end{itemize}
Notice that whenever $f$ is an $\mc{M}$-definable function, we can find $r \in \mathbf{B}_2$ and a large set $X$ such that \[\mc{M} \models \forall n, m \in X (f(n) = f(m) \Longleftrightarrow (n, m) \in \alpha_X(r)).\] 

We summarize the above as follows. Suppose $p(x) = \tp(\anglebracket{a, b} / \mc{M})$ and let $f$ be an $\mc{M}$-definable function. Then there is $\phi(x) \in \mc{L}(\mc{M})$ defining a ``large" set $X$ and $r \in \mathbf{B}_2$ such that $\phi(x) \in p(x)$ and \[\mc{M} \models \forall x, y [(\phi(x) \wedge \phi(y)) \Longrightarrow (f(x) = f(y) \Longleftrightarrow (x, y) \in \alpha_X(r))].\]
\end{ex}

Examples like this one provide the motivation for the following definitions. They were first introduced by Schmerl in \cite{schmerl_lat_86}, and they have been refined over the years.

\begin{defn} Let $L$ be a finite lattice, $X$ a set, and $\alpha : L \into \Eq(X)$ a representation.
\begin{enumerate}
    \item Let $Y \subseteq X$. Then $\alpha | Y : L \into \Eq(Y)$ is the pseudo-representation given by $(\alpha | Y)(r) = \alpha(r) \cap Y^2$ for each $r \in L$.
    \item Let $\beta : L \into \Eq(Y)$ be a pseudo-representation. Then $\alpha \cong \beta$ ($\alpha$ is \emph{isomorphic} to $\beta$) if there is a bijection $f : X \into Y$ such that for each $r \in L$, $(x, y) \in \alpha(r)$ if and only if $(f(x), f(y)) \in \beta(r)$.
    \item Let $\Theta \in \Eq(X)$. $\Theta$ is \emph{canonical} for $\alpha$ if there is $r \in L$ such that for all $x, y \in X$, $(x, y) \in \Theta$ if and only if $(x, y) \in \alpha(r)$.
    \item $\alpha$ has the \emph{$0$-canonical partition property}, or is \emph{$0$-CPP}, if for each $r \in L$, $\alpha(r)$ does not have exactly two classes.
    \item $\alpha$ is \emph{$(n+1)$-CPP} if, for each $\Theta \in \Eq(X)$ there is $Y \subseteq X$ such that $\alpha | Y$ is an $n$-CPP representation and $\Theta \cap Y^2$ is canonical for $\alpha | Y$.
\end{enumerate}
\end{defn}

Using these definitions, let us examine Example \ref{b2} once more from the other direction. Let $\mc{M} \models \PA$ and $X = [M]^2 = \{ \anglebracket{x, y} : x < y \}$. Then the representation $\alpha : \mathbf{B}_2 \into \Eq(X)$ given in Example \ref{b2} is $n$-CPP for each $n \in \omega$. To see this, first recall the Canonical Ramsey Theorem for pairs ($\mathsf{CRT}^2$): for every $f : [\omega]^2 \into \omega$, there is an infinite $X \subseteq \omega$ such that $f$ is canonical on $[X]^2$. That is, one of the following holds:
\begin{itemize}
    \item $f$ is one to one on $[X]^2$,
    \item $f$ is constant on $[X]^2$,
    \item for all  $\anglebracket{x_1, y_1}, \anglebracket{x_2, y_2} \in [X]^2$, $f(x_1, y_1) = f(x_2, y_2)$ if and only if $x_1 = x_2$, or
    \item for all $\anglebracket{x_1, y_1}, \anglebracket{x_2, y_2} \in [X]^2$, $f(x_1, y_1) = f(x_2, y_2)$ if and only if $y_1 = y_2$.
\end{itemize} This result is  due to Erd\H{o}s and Rado \cite{erdos-rado}, and is a consequence of Ramsey's Theorem for 4-tuples. One can formalize this result in $\PA^*$, so that if $\mc{M} \models \PA$ and $f : [M]^2 \into M$ is $\mc{M}$-definable, there is an $\mc{M}$-definable, unbounded $Y$ such that $f$ is canonical on $[Y]^2$. Notice, then, that for such a set $Y$, $\alpha | [Y]^2 \cong \alpha$.

Clearly, $\alpha$ is $0$-CPP. Moreover, if $\alpha$ is $n$-CPP, then by $\mathsf{CRT}^2$, for each $\Theta \in \Eq(X)$, there is $Y \subseteq X$ such that $\alpha \cong \alpha | Y$ and $\Theta \cap Y^2$ is canonical for $\alpha | Y$. Since $\alpha | Y \cong \alpha$, then $\alpha | Y$ is $n$-CPP, and therefore $\alpha$ is $(n+1)$-CPP.

Given this $\alpha$, one can construct an elementary extension $\mc{N}$ of $\mc{M}$ such that $\Lt(\mc{N} / \mc{M}) \cong \mathbf{B}_2$. The idea is to construct a type $p(x)$ ensuring that, for each $\mc{M}$-definable function $f$, there is some definable $Y \subseteq X$ such that the equivalence relation induced by $f$ is canonical for $\alpha$ on $Y$, and that the defining formula for $Y$ is in $p(x)$.

To construct this type, we will construct an infinite descending sequence of ``large" sets. Let $X_0 = [M]^2$. Enumerate the $\mc{M}$-definable equivalence relations $\Theta_0, \Theta_1, \ldots$. Given $\alpha | X_i$ and $\Theta_i$, we use $\mathsf{CRT}^2$ to find $X_{i+1} \subseteq X_i$  such that $\alpha | X_{i_1} \cong \alpha | X_i$ and $\Theta_i$ is canonical for $\alpha | X_{i+1}$. Finally,  we let $p(x)$ be the type $$\{ \phi(x) \in \mc{L}(\mc{M}) : \text{there is } i \in \omega \text{ such that } \mc{M} \models \forall x (x \in X_i \rightarrow \phi(x)) \}.$$

We show that that $p(x)$ is a complete type. This is, essentially, due to the fact that each $\alpha | X_i$ is 0-CPP. That is, given $\phi(x) \in \mc{L}(\mc{M})$, consider the equivalence relation $\Theta$ given by $(x, y) \in \Theta$ iff $\mc{M} \models \phi(x) \Longleftrightarrow \phi(y)$. Let $\Theta = \Theta_i$, and notice that since $\Theta$ is canonical for $\alpha | X_{i+1}$  and $\Theta$ has at most two equivalence classes, it must be the case that $\Theta \cap X_{i+1}^2$ is trivial.

Let $c$ realize $p(x)$. We show why $\Lt(\mc{M}(c) / \mc{M}) \cong \mathbf{B}_2$. Because the pairing function is one to one, there are $a$ and $b$ such that $\mc{M}(c) \models c = \anglebracket{a, b}$. We use the same names as in the lattice $\mathbf{B}_2$ (see Figure \ref{fig1}) suggestively. For each $d \in \mc{M}(c)$, let $f$ be an $\mc{M}$-definable function such that $\mc{M}(c) \models f(c) = d$, and let $\Theta$ be the equivalence relation induced by $f$. Then there is a simple argument that exactly one of the following must hold:
\begin{itemize}
    \item $\mc{M}(d) = \mc{M}$,
    \item $\mc{M}(d) = \mc{M}(a)$,
    \item $\mc{M}(d) = \mc{M}(b)$, or
    \item $\mc{M}(d) = \mc{M}(c)$.
\end{itemize} This is proved case by case by finding $r$ and $i$ such that $\Theta \cap X_i^2 = \alpha(r) \cap X_i^2$ (by canonicity). Moreover, it is clear that $\mc{M} \prec  \mc{M}(a), \mc{M}(b) \prec \mc{M}(c)$. One checks that $\mc{M}(a) \cap \mc{M}(b) = \mc{M}$.

The definitions used above relativize to a model $\mc{M} \models \PA$. Suppose $X \subseteq M$ is $\mc{M}$-definable and $\alpha : L \into \Eq(X)$ is a representation. Then $\alpha$ is an \emph{$\mc{M}$-representation} if $\alpha$ is $\mc{M}$-definable. If $X \in \Def(\mc{M})$, then by $\Eq^{\mc{M}}(X)$ we mean the lattice of $\mc{M}$-definable equivalence relations on $X$. Similarly, the notion of $n$-CPP representations formalizes in $\lpa$ as well; in such cases, one only considers representations over $\mc{M}$-finite sets, and we quantify over the equivalence relations $\Theta \in \Eq^{\mc{M}}(X)$. That is, there is an $\lpa$-formula $cpp_L(x)$ asserting that $L$ has an ($\mc{M}$-finite) $x$-CPP representation.

In the construction of the type $p(x)$ above, the two important ingredients needed at each step are:
\begin{itemize}
    \item ensure that each $\alpha | X_i$ is $0$-CPP, and
    \item ensure that each $\Theta_i$ is canonical for $\alpha | X_{i+1}$.
\end{itemize} This observation naturally leads to the following definitions and results in \cite{schmerl_n5} by James Schmerl, further refining the technique.

\begin{defn}[{\cite[Definition 1.3]{schmerl_n5}}] Let $\mc{M} \models \PA$ and $L$ a finite lattice. $\mc{C}$ is an \emph{$\mc{M}$-correct set of representations of $L$} if each $\mc{C}$ is a nonempty set of $0$-CPP $\mc{M}$-representations of $L$ and whenever $\alpha : L \into \Eq(X) \in \mc{C}$ and $\Theta \in \Eq^{\mc{M}}(X)$, there is $Y \subseteq X$ such that $\alpha | Y \in \mc{C}$ and $\Theta \cap Y^2$ is canonical for $\alpha | Y$.
\end{defn}

Returning to Example \ref{b2}, we notice that if $\mc{C}$ is the collection of $\alpha | Y$ such that $Y$ is infinite, $\mc{M}$-definable and $\alpha | Y \cong \alpha$, then $\mc{C}$ is an $\mc{M}$-correct set of representations of $\mathbf{B}_2$. Additionally, one observes that $\PA \vdash cpp_{\mathbf{B}_2}(n)$ for each $n \in \omega$. If $\mc{M} \models \PA$ is nonstandard, by overspill there is a nonstandard $c$ such that $\mc{M} \models cpp_{\mathbf{B}_2}(c)$. Then the collection of all $\mc{M}$-representations of $\mathbf{B}_2$ that are $x$-CPP for some nonstandard $x$ is also $\mc{M}$-correct.

\begin{thm}[{\cite[Theorem 1.4]{schmerl_n5}}]
Let $\mc{M} \models \PA$ and $L$ be a finite lattice. Then:
\begin{enumerate}
    \item If there is $\mc{N}$ such that $\mc{M} \prec \mc{N}$ and $\Lt(\mc{N} / \mc{M}) \cong L$, then there is an $\mc{M}$-correct set of representations of $L$.
    \item If $\mc{M}$ is countable and there is an $\mc{M}$-correct set of representations of $L$, then there is $\mc{N} \succ \mc{M}$ such that $\Lt(\mc{N} / \mc{M}) \cong L$.
\end{enumerate}
\end{thm}

\section{Representations of Ranked Lattices}\label{ranked}

In this section, we extend the definition of $\mc{M}$-correct sets of representations of a lattice to ranked lattices as in Section \ref{ranked-def}. Let us first consider Example \ref{b2} in the context of ranked lattices. 

Let $\mc{M} \prec \mc{N}$ be such that $\Lt(\mc{N} / \mc{M}) \cong \mathbf{B}_2$, and let $a, b \in N$ be such that $\Lt(\mc{N} / \mc{M}) = \{ \mc{M}, \mc{M}(a), \mc{M}(b), \mc{M}(a, b) = \mc{N} \}$ (as we did previously, we use $a$ and $b$ suggestively to correspond with $a, b \in \b2$). Because the rankset is linearly ordered, it must be the case that either $\mc{M}(a) \prec_{\text{cof}} \mc{N}$ or $\mc{M}(b) \prec_\text{cof} \mc{N}$. Without loss of generality, assume that $\mc{M}(b) \prec_\text{cof} \mc{N}$. Therefore, the possible ranksets of $\Lt(\mc{N} / \mc{M})$ are:
\begin{itemize}
    \item $ \{ \mc{N} \}$ (if $\mc{M} \prec_\text{cof} \mc{N}$),
    \item $\{ \mc{M}, \mc{M}(a), \mc{N} \}$ (if $\mc{M} \prec_\text{end} \mc{M}(a) \prec_\text{end} \mc{N}$), or,
    \item $\{ \mc{M}(a), \mc{N} \}$ (if $\mc{M} \prec_\text{cof} \mc{M}(a) \prec_\text{end} \mc{N}$). 
\end{itemize}

To study an example of a mixed extension, we consider the case where $\mc{M} \prec_\text{cof} \mc{M}(a) \prec_\text{end} \mc{N}$, and $\mc{M}(b) \prec_\text{cof} \mc{N}$. Recall that we referred to a set $X$ as \emph{large} if it is infinite, $\mc{M}$-definable, and $\anglebracket{a, b} \in X^\mc{N}$. We look for properties to motivate a notion of ``largeness" in a mixed extension.

Let $X = \{ \langle x, y \rangle : x < y \}$ and $\alpha_X$ be as defined in Example \ref{b2}. That is, $\alpha_X : \b2 \to \Eq(A)$ is defined so that $\alpha_X(a)$ is the equivalence relation induced by $\pi_1$ (the projection onto the first coordinate, i.e., $\anglebracket{x, y} \mapsto x$), and $\alpha_X(b)$ is the equivalence relation induced by $\pi_2$ (projection onto the second coordinate). Then there is an infinite, $\mc{M}$-definable $Y \subseteq A$ such that, letting $\alpha_Y$ be $\alpha_X | Y$:

\begin{enumerate}
    \item\label{cof0} There is a bounded set of representatives of the collection of all $\alpha_Y(a)$-classes.
    \item\label{enda}There is an unbounded $\alpha_Y(a)$-class (that is, a class which contains unboundedly many $\alpha_Y(1)$-classes).
    \item\label{cofb} Every $\alpha_Y(b)$-class is $\mc{M}$-finite. 
\end{enumerate}

To see these, notice that since $\mc{M} \prec_\text{cof} \mc{M}(a)$, then there is $m \in M$ such that $\mc{M} \models a < m$. Let $Y = \{ \langle x, y \rangle: x < m$ and $x < y \}$. The following statements are easily verified:

\begin{itemize}
    \item There are (exactly) $m$ $\alpha_Y(a)$-classes (one for each $x < m$),
    \item each $\alpha_Y(a)$-class is unbounded, and,
    \item each $\alpha_Y(b)$-class has at most $m$ elements.
\end{itemize} Moreover, since $\mc{N} \models a < m$ and $a < b$, then $\mc{N} \models \anglebracket{a, b} \in Y^{\mc{N}}$. Similarly to Example \ref{b2}, if we define a set $Y$ to be large if it is $\mc{M}$-definable, satisfies properties $\eqref{cof0}-\eqref{cofb}$ and $\mc{N} \models \anglebracket{a, b} \in Y^{\mc{N}}$, then whenever $Y$ is large and $\Theta \in \Eq^\mc{M}(Y)$, there is large $Z \subseteq Y$ such that $\Theta$ is canonical for $\alpha_Z$.  

We point out here some specific features of this representation that follow from properties (1)-(3). These features, it turns out, need to be present for any representation of the ranked lattice $(\b2, \rho)$, where the rankset of $\rho$ is $\{ a, 1 \}$. First, $0 < a = \rho(0)$. In every representation, $\alpha(0)$ is trivial. In this case, notice that, in a sense, the lone $\alpha_Y(0)$-class ($Y$) splits into boundedly many $\alpha_Y(a)$-classes; we will make this notion more precise in the below definition. It is easy to see that each $\alpha_Y(a)$-class is of the form $\{ \anglebracket{x, y} : x < y \}$ for some fixed $x \in M$ such that $\mc{M} \models x < m$ and that  $Y$ is the union of all of these classes.

Secondly, $0 < b$ but $b \not< \rho(0)$. Notice now that $Y$ is not the union of an $\mc{M}$-bounded set of $\alpha_Y(b)$-classes. Again, we can see this because each $\alpha_Y(b)$-class is of the form $\{ \anglebracket{x, y} : x < y \}$ for some fixed $y \in M$.

Lastly, $b < 1 = \rho(b)$; again, see that each $\alpha_Y(b)$-class splits into boundedly many $\alpha_Y(1)$-classes; in other words, each $\alpha_Y(b)$-class is $\mc{M}$-finite.

\begin{defn}[{\cite[Definition 1.6]{schmerl_n5}}]Let $\mc{M} \models \PA$ and $(L, \rho)$ a finite ranked lattice.
\begin{enumerate}
    \item If $A \in \Def(\mc{M})$ and $\Theta \in \Eq(A)$ is $\mc{M}$-definable, a set $\mc{E}$ of $\Theta$ classes is $\mc{M}$-bounded if there is a bounded $I \in \Def(\mc{M})$ such that $I \cap X \neq \emptyset$ for each $X \in \mc{E}$.
    \item $\alpha : L \to \Eq(A)$ is an \emph{$\mc{M}$-representation of $(L, \rho)$} if $\alpha$ is an $\mc{M}$-representation of $L$ (that is, $\alpha$ is $\mc{M}$-definable) and whenever $r \leq s \in L$, $s \leq \rho(r)$ if and only if every $\alpha(r)$-class is the union of an $\mc{M}$-bounded set of $\alpha(s)$-classes.
    \item $\mc{C}$ is an \emph{$\mc{M}$-correct set of representations of $(L, \rho)$} if $\mc{C}$ is an $\mc{M}$-correct set of representations of $L$ and each $\alpha \in \mc{C}$ is an $\mc{M}$-correct representation of $(L, \rho)$.
\end{enumerate}
\end{defn}

Notice that in the above example,  whenever $X \subseteq A$ is large, $\alpha_X$ is an $\mc{M}$-representation of $(\b2, \rho)$, where $\rho$ is the ranking whose rankset is $\{ a, 1 \}$.

\begin{thm}[{\cite[Theorem 1.7]{schmerl_n5}}]Suppose $\mc{M} \models \PA$ and $(L, \rho)$ is a finite ranked lattice.
\begin{enumerate}
    \item If there is $\mc{N}$ such that $\mc{M} \prec \mc{N}$ and $\Ltr(\mc{N} / \mc{M}) \cong (L, \rho)$, then there is an $\mc{M}$-correct set of representations of $(L, \rho)$.
    \item If $\mc{M}$ is countable and there is an $\mc{M}$-correct set of representations of $(L, \rho)$, then there is $\mc{N} \succ \mc{M}$ such that $\Ltr(\mc{N} / \mc{M}) \cong (L, \rho)$.
\end{enumerate}
\end{thm}

We turn now to an example of an $\mc{M}$-correct set of representations of $\m3$ (see Figure \ref{fig1}). First we define a representation $\alpha : \m3 \to \Eq(3)$, where $3$ is the set $\{ 0, 1, 2 \}$. Define this representation as follows:

\begin{itemize}
    \item the equivalence classes of $\alpha(a)$ are $\{0 \}$ and $\{ 1, 2 \}$,
    \item the equivalence classes of $\alpha(b)$ are $\{ 0, 2 \}$ and $\{ 1 \}$, and,
    \item the equivalence classes of $\alpha(c)$ are $\{ 0, 1 \}$ and $\{ 2 \}$.
\end{itemize}

Let $\mc{M}$ be a countable, nonstandard model and $m \in M$ (standard or nonstandard). Let $3^m$ refer to the set of (codes of) $\mc{M}$-finite sequences $s$ whose length is $m$ and, for each $i < m$, $(s)_i \in \{ 0, 1, 2 \}$. Define $\alpha^m : \m3 \to \Eq(3^m)$ by letting $(s, t) \in \alpha^m(r)$ if and only if $((s)_i, (t)_i) \in \alpha(r)$ for each $i < m$. Notice that each of these representations is a representation of $(\m3, \rho)$ where $\rho(r) = 1$ for each $r \in \m3$ (that is, there are $\mc{M}$-boundedly many $\alpha^m(a)$ classes, each one is $\mc{M}$-bounded, etc).

It turns out that the set of $\alpha^m$ when $m$ is nonstandard forms an $\mc{M}$-correct set of representations of $(\m3, \rho)$. This is not obvious: one appeals to a generalization of the Hales-Jewett Theorem due to Pr\"{o}mel and Voigt. In fact, this phenomenon can be generalized to any finite lattice $L$ which can be represented as a congruence lattice of a finite algebra.

\subsection{Congruence Lattices}

An \emph{algebra} is a structure of the form $(A, \anglebracket{f_i : i \in I})$, where $A$ is a set, $I$ is an index set, and for each $i \in I$, there is some $n \in \omega$ such that $f_i : A^n \to A$ (we allow for $n$ to be $0$, in which case such an $f_i$ is a constant).

If $\mc{A} = (A, \anglebracket{f_i : i \in I})$ is an algebra, then a \emph{congruence} is an equivalence relation $\theta$ on $A$ which commutes with all of the $f_i$. That is, for each $i \in I$, if $f_i : A^n \to A$ and $\bar{a}, \bar{b}$ are tuples of length $n$ such that $(a_j, b_j) \in \theta$ for all $j < n$, then $(f_i(\bar{a}), f_i(\bar{b})) \in \theta$. The set of all congruences on an algebra $\mc{A}$ is denoted $\Cg(\mc{A})$, and forms a sublattice of $\Eq(A)$.

In the definition below, for a lattice $L$, we let $L^d$ be its dual, i.e., $L$ with its ordering reversed.

\begin{defn}{\cite[Definition 4.5.7]{ks}} Let $L$ be a finite lattice and $\alpha : L \to \Eq(A)$ a representation. $\alpha$ is a \emph{congruence representation} if there is an algebra $\mc{A}$ such that $\alpha$ is an isomorphism of $L$ and $\Cg(\mc{A})^d$. 
\end{defn}

Every algebraic lattice is isomorphic to a congruence algebra \cite{gratzer-schmidt}. It is a well-known open question in universal algebra whether every finite lattice has a \emph{finite} congruence representation; that is, if it can be represented as $\Cg(\mc{A})$ for a finite algebra $\mc{A}$. This problem is referred to as the \emph{finite lattice representation problem}. In conjunction with the next result due to Schmerl (\cite{schmerl_finite_lat93}), a positive result to the finite lattice representation problem implies a positive result for the restriction of the lattice problem for models of $\PA$ to finite lattices.

\begin{thm}
Let $L$ be a finite lattice which has a finite congruence representation. Then every countable nonstandard $\mc{M} \models \PA$ has a cofinal elementary extension $\mc{N}$ such that $\Lt(\mc{N} / \mc{M}) \cong L$.
\end{thm}

As mentioned before, the proof of this is a generalization of the result for $\m3$ given above. For details, see \cite[Section 4.5]{ks}, and, in particular, \cite[Theorem 4.5.27 and Corollary 4.5.28]{ks}.

\section{Diversity} 
The early results of Gaifman and Mills mentioned in Section 3 relied heavily on applications of minimal and end-extensional types. For a given $\mc{M}$, if $\mc{N}$ is generated over $\mc{M}$ by a set $A$ of elements realizing the same minimal \elm-type, then \Lat\ is isomorphic to the Boolean algebra of subsets of $A$ and the isomorphism type over $\mc{M}$ of each model in \Lat\ is determined by the cardinality of its set of generators. In particular, if $A$ is finite of cardinality $n$, then there are exactly $n+1$ isomorphism types of models in \Lat. If $\mc{N}$ is generated over $\mc{M}$ by a finite set of elements realizing mutually independent minimal types, then no pair of distinct models in \Lat\ are isomorphic.

In \cite{schmerl_div1} Schmerl asked if there is more that can be said about the diversity of isomorphism types in \Lat. He called a model  \emph{diverse} if no two of its elementary submodels are isomorphic, and called an extension $\mc{N}$ of $\mc{M}$ \emph{diverse} if no  two models in \Lat\ are isomorphic over $\mc{M}$. 

In \cite{schmerl_div1}, Schmerl observes that $\mc{M}$ is diverse if and only if no two distinct elements have the same type. For such structures (for any language) Ali Enayat coined the name \emph{Leibnizian}. His article \cite{leibnizian} is devoted to Leibnizian models of set theory.

All results in the rest of this section are from \cite{schmerl_div1}. Schmerl notes that the first theorem suggests that constructing models that are not diverse is more difficult than constructing diverse ones. 

\begin{thm}
\begin{enumerate}
\item If $\mc{M}$ is not an elementary extension of the standard model and \lat\ is finite, then there is a diverse $\mc{N}$ such that $\mc{M}\equiv\mc{N}$ and $\Lt(\mc{M})\cong \Lt(\mc{N})$.
\item Let $L$ be a finite lattice. If $\mc{M}$ is nonstandard and has an elementary extension $\mc{N}$ such that \Lat$\cong L$, then $\mc{M}$ has a cofinal diverse  extension $\mc{N}$ such that \Lat$\cong L$.
\end{enumerate}
\end{thm}
Because a model is diverse just in case it is a diverse extension of its prime elementary submodel, (1) above is an easy consequence of the second.  The proof of (2) involves the full power of the CPP-representations theory,  invoking at one point a canonical partition theorem of Pr\"omel and Voigt (see \cite[page 118]{ks} for the statement of the theorem).

If $E$ is an equivalence relation on a lattice $L$, then  $(L,E)$ is called an \emph{equivalenced lattice}. By \latp\ we denote $(\Lt(\mc{N}), E)$, where $E$ is the isomorphism relation, and \Latp\ is $(\Lt(\mc{N} / \mc{M}), E)$, where $E$ is the  isomorphism relation for the isomorphisms that fix $M$ pointwise. With this notation, the basic lattice problem gets generalized to
\begin{prob} For which finite equivalenced lattices $(L,E)$ are there models $\mc{M}$ such that  \latp$\cong (L,E)$ and for which finite equivalenced lattices $L$ are there models $\mc{M}$ and $\mc{N}$  such that \Latp$\cong (L,E)$?
\end{prob}
Right from the start one can see that some equivalenced lattices cannot be represented as substructure and interstructure lattices for obvious reasons. For example if $L$ is $\n5$ and $(b,c)\in E$ (Figure 2), then, as Schmerl writes: ``it would not be at all reasonable to expect such models." This leads to the definition of reasonable equivalenced lattice. 

Before stating the definition, let us observe that each model of the form $\mc{M}(a)$  is rigid over $\mc{M}$. This is a consequence of Ehrenfeucht's lemma \cite{ehrenfeucht} (see \cite[Theorem 1.7.2]{ks}), which says that  for all $b\in \mc{M}(a)$, if $\tp(a/M)=\tp(b/M)$ then $a=b$. It follows that if $\mc{M}(a_1)$ and $\mc{M}(a_2)$ in \Lat\ are isomorphic over $\mc{M}$ and $\tp(a_1/M)=\tp(a_2/M)$, then there is a unique isomorphism $F:\mc{M}(a_1)\into \mc{M}(a_2)$ such that $F(a_1)=a_2$.

If $(L,E)$ is a finite equivalenced lattice and $F:L\into$ \Latp\ is an isomorphism, then we can define a linear order $\lhd$ of $L$ as follows. For each $r\in L$ select a generator $a_r$ of $F(r)$ over $\mc{M}$ so that for all $r$ and $s$, if $(r,s)\in E$, then $\tp(a_r/M)=\tp(a_s/M)$ and define $r\lhd s$ iff $a_r<a_s$.

Schmerl calls an equivalenced lattice $(L,E)$ \emph{reasonable} if there is a linear ordering $\lhd$  on $L$ such that whenever $(a,b)\in E$ and $I$ and $J$ are the principal ideals of $L$ generated by $a$ and $b$ respectively, then there is an isomorphism $f:(I,\lhd\cap I^2)\into (J,\lhd\cap J^2)$ such that $(x,f(x))\in E$ for all $x\in I$. One can directly check that if \Lat\ is finite,  then  \Latp\ is reasonable. 

The definition given above is from the later Schmerl's paper \cite{schmerl_div2}. Under a less restrictive definition given in \cite{schmerl_div1}, Schmerl proves that for every reasonable equivalenced Boolean lattice $(\bn,E)$, every nonstandard model $\mc{M}$ has a cofinal extension such that \Latp\ is isomorphic to $(\bn,E)$. In \cite{schmerl_div2} this is generalized to: under the revised definition given above, for every reasonable equivalenced distributive lattice $(L,E)$, every nonstandard model $\mc{M}$ has a cofinal extension $\mc{N}$ such that \Latp\ is isomorphic to $(L,E)$. 

 Among many open problems about diversity, this is probably the simplest:
 \begin{prob} Let $\mc{M}$ be countable and nonstandard. Is there a nondiverse extension $\mc{N}$ such that either \Lat\ $=\n5$ or \Lat\ $= \m3$?
 \end{prob}
 
\section{Coded sets and distributive lattices}

Recall that if $\mathcal{N}$ is an end extension of $\mathcal{M}$, then  $\Cod(\mathcal{N} / \mathcal{M})$ is $\{ X \cap M : X \in \Def(\mathcal{N}) \}$. We referred to this as the family of \emph{coded sets} in the extension. Schmerl \cite{schmerlmin} characterized exactly which families of subsets of a model of $\PA$ (of arbitrary cardinality) can appear as $\Cod(\mathcal{N} / \mathcal{M})$ when $\mathcal{N}$ is a minimal elementary end extension of $\mathcal{M}$; i.e., an extension where $\Lt(\mathcal{N} / \mathcal{M}) \cong \mathbf{2}$. 

\begin{thm}[{\cite[Theorem 3]{schmerlmin}}]
If $\mc{M} \models \PA$ and $\mathfrak{X} \subseteq \mc{P}(\mc{M})$, the following are equivalent:
\begin{enumerate}
    \item There is a countably generated extension $\mc{N} \succ_\text{end} \mc{M}$ such that $\Cod(\mc{N} / \mc{M}) = \mathfrak{X}$ and every set that is $\Pi^0_1$-definable in $(\mc{M}, \mathfrak{X})$ is the union of countably many $\Sigma^0_1$-definable sets.
    \item There is a minimal extension $\mc{N} \succ_\text{end} \mc{M}$ such that $\Cod(\mc{N} / \mc{M}) = \mathfrak{X}$.
\end{enumerate}
\end{thm}

Previously, Schmerl \cite{schmerlcoded} characterized the families of sets that can appear as $\Cod(\mc{N} / \mc{M})$ in any countably generated elementary end extension $\mc{N}$ of a model $\mc{M}$, so this result completes the picture for minimal extensions.

A simple construction\footnote{Schmerl (2017) via private communication} shows that the same characterization of those coded sets holds for any finite distributive lattice $D$. That is:

\begin{prop} Let $\mathcal{M} \models \PA$ and $\mathfrak{X} \subseteq \mathcal{P}(\mathcal{M})$. The following are equivalent:
\begin{enumerate}
    \item\label{min} There is a minimal elementary end extension $\mathcal{N}$ of $\mathcal{M}$ such that $\Cod(\mathcal{N} / \mathcal{M}) = \mathfrak{X}$.
    \item\label{fdl} For any finite distributive lattice $D$, there is $\mathcal{N} \succ_\text{end} \mathcal{M}$ such that $\Lt(\mathcal{N} / \mathcal{M}) \cong D$ and $\Cod(\mathcal{N} / \mathcal{M}) = \mathfrak{X}$.
\end{enumerate}
\end{prop}

Before we describe Schmerl's proof of this result, we need the following facts about finite distributive lattices (see \cite[Section 4.3]{ks}). Let $L$ be a lattice and $a \in L$. The \emph{$a$-doubling extension} of $L$ is the sublattice $L^\prime$ of $L \times \mathbf{2}$ (ordered lexicographically; i.e., $(r, i) \leq (s, j)$ iff $r \leq s$ and $i \leq j$) defined as $\{ (r, i) \in L \times \mathbf{2} : i = 0$ or $r \geq a \}$. \cite[Theorem 4.3.6]{ks} states that a finite lattice $L$ is distributive if and only if there is sequence $L_0, \ldots, L_n$ of lattices such that $L_0$ is the one-element lattice, $L_n \cong L$, and each $L_{i+1}$ is a doubling extension of $L_i$.

\begin{proof}
$\eqref{fdl} \implies \eqref{min}$ is clear since $\mathbf{2}$ is a finite distributive lattice, so assume that $\mc{M} \prec \mc{N}$ is a minimal elementary end extension and $D$ is a finite distributive lattice. Then let $L_0, \ldots, L_n \cong D$ be the finite sequence of doubling extensions as stated above.

We show that whenever $\mc{M}_1 \prec_\text{end} \mc{M}_2$ is such that $\Lt(\mc{M}_2 / \mc{M}_1)$ is finite, then for any $\mc{K} \in \Lt(\mc{M}_2 / \mc{M}_1)$, $\mc{M}_2$ has an elementary end extension $\mc{N}^\prime$ such that $\Lt(\mc{N}^\prime / \mc{M}_1)$ is isomorphic to the $\mc{K}$-doubling extension of $\Lt(\mc{M}_2 / \mc{M}_1)$. The conclusion to $\eqref{fdl}$ follows, since if $\mc{M}_1 \prec_\text{end} \mc{M}_2 \prec_\text{end} \mc{N}^\prime$, then $\Cod(\mc{N}^\prime / \mc{M}_1) = \Cod(\mc{M}_2 / \mc{M}_1)$.

To find such an $\mc{N}^\prime$, first notice that since $\Lt(\mc{M}_2 / \mc{M}_1)$ is finite, then $\mc{K} = \mc{M}_1(a)$ for some $a \in K$. Let $\mc{M}_0 = \Scl(a)$. By \cite[Theorem 4.3.2]{ks}, since $\mc{M}_0$ is countable, $\mc{M}_2$ has an elementary end extension $\mc{N}^\prime$ such that $\Lt(\mc{N}^\prime)$ is isomorphic to the $\mc{M}_0$-doubling extension of $\Lt(\mc{M}_2)$. This $\mc{N}^\prime$ is as required; that is, $\Lt(\mc{N}^\prime / \mc{M}_1)$ is isomorphic to the $\mc{K}$-doubling extension of $\Lt(\mc{M}_2 / \mc{M}_1)$. This follows immediately because $\mc{M}_0 \vee \mc{M}_1 = \mc{K}$.
\end{proof}

Note that if $\mathcal{M}$ is countable and $\mathfrak{X} = \Def(\mc{M})$, the conclusion to (1) holds since every countable model has conservative minimal elementary end extensions. Since the pentagon lattice $\mathbf{N}_5$ cannot be realized as the interstructure lattice of a conservative end extension, we do not have the same characterization for non-distributive lattices.
 
\section{Countable \rs\ models}
 Up to this point, this survey was about the  the problem of finding, for a given lattice $L$, a substructure or interstructure lattice representation of $L$ with or without some additional properties. In this section we will briefly discuss a dual problem for a particular class of models. Given a \ct\ \rs\ model, what can we say about its substructure lattice? 

A model $\mc{M}$ is \rs\ if it is saturated with the respect to computable types with finite numbers of parameters. The \emph{standard system} of a model $\mc{M}$, $\SSy(\mc{M})$, is the set of standard parts of the definable subsets of $\mc{M}$, i.e., $\SSy(\mc{M})=\Cod(\mc{M}/\mathbb N)$, where $\mathbb N$ is the standard model. It is not difficult to prove that any two  \ct\ \rs\ models are isomorphic if and only if they are elementarily equivalent and they have the same standard system. 

{\bf For the rest of this section, let $\mc{M}$ be \ct\ and \rs.} 

The lattice \lat\ is immense. If the domain of a model $\mc{K}$ in \lat\ is an initial segment of $M$ we call $\mc{K}$ an \emph{elementary cut}. Henryk Kotlarski proved that  the set of elementary cuts of $\mc{M}$ ordered by inclusion is isomorphic to $2^\omega$ with the lexicographic ordering.  Moreover, he proved that the set of elementary cuts that are \rs\ is uncountable, dense, and is closed under infinite unions but not under infinite intersections \cite{kot84}. All  \rs\ elementary cuts have the same standard system; hence, they  are isomorphic to one another. Yet, there are continuum many first-order theories  of pairs $(\mc{M}, K)$, where $\mc{K}$ is a \rs\ elementary cut \cite{smo82}.

An elementary cut of $\mc{M}$ is not \rs\ if and only if it is a closure under initial segment of the Skolem closure of single element of $M$. Every \rs\ model realizes countably many mutually independent minimal types. It follows that there are countably many isomorphism types of elementary cuts that are not \rs. Moreover, each cut that is not \rs, except for the closure under initial segment of  $\scl(0)$, has countably  many automorphic images; hence the isomorphism relation of \lat\ restricted to such  cuts has countably many countable equivalence classes.

Elementary cuts are linearly ordered by inclusion; hence they form a distributive sublattice of \lat. To show that \lat\ is not distributive we will use some facts about  cofinal extensions.  

In \cite{schmerl_almost_min}, Schmerl calls the extension $\mc{M}\prec \mc{N}$ \emph{almost minimal} if $\mc{N}$ is an end extension and for every $\mc{K}$, if  $\mc{M}\prec\mc{K}\prec\mc{N}$, then $\mc{K}\cof\mc{N}$. For a lattice $L$, ${\mathbf 2}\oplus L$ is the extension of $L$ obtained by adding with a new 0 element appended below ${\mathbf 0}_L$. By the Blass condition, if $\mc{N}$ is an almost minimal extension of $\mc{K}$ and $\Lt(\mc{N}/\mc{K})$ is finite, then $\Lt(\mc{N}/\mc{L})\cong {\mathbf 2}\oplus L$ for some finite lattice $L$. The following is an abridged version of \cite[Corollary 3.3]{schmerl_almost_min}.

\begin{thm}\label{notdist}
Let $\mc{K}$ be a \ct\ model of $\PA$ and let $L$ be a finite lattice. The following are equivalent.
\begin{enumerate}
    \item $\mc{K}$ has an elementary extension $\mc{N}$ such that $\Lt(\mc{N}/\mc{K})\cong {\mathbf 2}\oplus L$.
  \item $\mc{K}$ has an elementary almost minimal extension $\mc{N}$ such that $\Lt(\mc{N}/\mc{K})\cong {\mathbf 2}\oplus L$.  \end{enumerate}
\end{thm}
For our discussion here it is important that if $\mc{N}$ is an almost minimal extension of a nonstandard $\mc{K}$ then $\SSy(\mc{K})=\SSy(\mc{N})$. In particular, if  $\mc{K}$ and $\mc{N}$ are \ct\ and  \rs, then $\mc{K}\cong\mc{N}$,  which shows that ${\mathbf 2}\oplus L$ in (2) embeds into $\Lt(\mc{K})$.

We are going back to our countable \rs\ model $\mc{M}$.  Let $S$ be an inductive partial satisfaction class on $\mc{M}$ (see \cite[Definition 1.9.1]{ks}) and let $\mc{M}^*$ be $(\mc{M},S)$. Let $L$ be a finite lattice satisfying (1) in  the $\PA^*$ version of Theorem \ref{notdist}, i.e., for $\mc{K}=\mc{M}^*$, and let $\mc{N}^*$ be an almost minimal extension of $\mc{M}^*$ given by (2). By the remarks following the theorem, we get that ${\mathbf 2}\oplus L$ embeds into $\Lt(\mc{N}^*)$, hence it also embeds into $\Lt(\mc{N})$.
Applying this to $L=\m3$, we get the following corollary.

\begin{cor}\label{notdist_cor}
Lattices of elementary substructures of \rs\ models of $\PA$ are not distributive.
\end{cor}
It seems that the argument given above to prove Corollary \ref{notdist_cor} is an overkill. More directly, one can construct a type $p(x) \in \SSy(\mc{M})$, formalizing \cite[Theorem 4.5.21]{ks}, such that if $a \in M$ realizes $p(x)$, then $\Lt(\Scl(a)) \cong \mathbf{2} \oplus \m3$.

Much more can be said about cofinal submodels of $\mc{M}$. Here is a sample. It follows from \cite[Theorem 7.1]{ks12} that there is a set $\mathcal C$ of continuum many cofinal submodels of $\mc{M}$ such that for all $\mc{K}\in \mathcal C$:

\begin{itemize}
\item  $\mc{K}$ is isomorphic to $\mc{M}$; 
\item $\Lt(\mc{M}/\mc{K})$ is the three-element lattice $\mathbf 3$; 
\item  for each pair of distinct $\mc{K}_1$, $\mc{K}_2$ in $\mathcal C$,  $\Th(\mc{M},K_1)\not=\Th(\mc{M},K_2)$. 
\end{itemize}
In the recent paper \cite{schmerl_min_cof}, after much work, Schmerl improves this by replacing $\mathbf 3$ by $\mathbf 2$, i.e., for each $\mc{K}$, $\mc{M}$ is a minimal extension of $\mc{K}$.

At the end of \cite{kot84}, Kotlarski posed a general problem to describe the \st\  of \lat\ for \ct\ \rs\ $\mc{M}$. In particular, he asked whether \lat\ depends on $\mc{M}$. 

A partial answer to Kotlarski's question is given in \cite[Section 5]{ks95}. It solves the problem for \as\ models. A \rs\ model of $\PA$ is \emph{\as} if its standard system is closed under \ar\ comprehension. It is shown in \cite{ks95} that if $\mc{M}$ and $\mc{N}$ are \ct\ \as\ models of the same completion of $\PA$, then $\mc{M}\cong\mc{N}$ if and only if $\Lt(\mc{M})\cong\Lt(\mc{N})$. For the proof, Schmerl introduced a family of countably infinite, distributive lattices $D(X)$, one for each set of natural numbers $X$, such that:
\begin{enumerate}
    \item if $X$ and $Y$ are distinct then $D(X)$ and $D(Y)$ are not isomorphic;
    \item \label{lat-determines-ssy} if $\mc{M}$ is \as\  and there is $b\in M$ such that $D(X)\cong \Lt(\scl(b))$, then $X\in \SSy(\mc{M})$; 
    \item \label{ssy-determines-lat} if $\mc{M}$ is \rs\ and $X$ is in the standard system of $\mc{M}$, then there is $b\in M$ such that $D(X)\cong \Lt(\scl(b))$.
\end{enumerate} 
This proves the following theorem.
\begin{thm}\label{varlt} Let $T$ be a completion of $\PA$ and let $\mc{M}$ and $\mc{N}$ be \ct\ \as\ models of $T$. Then $\Lt(\mc{M})\cong\Lt(\mc{N})$ if and only if $\mc{M}\cong\mc{N}$.
\end{thm}

For \as\ models,  \eqref{lat-determines-ssy} and \eqref{ssy-determines-lat} above could be stated as an equivalence. The reason for their separation is that while the proof of \eqref{lat-determines-ssy} is relatively easy, it seems to require the full strength of arithmetic saturation. The fact \eqref{ssy-determines-lat} above holds for all \rs\ models, but its proof heavily depends on the methods developed by Schmerl in \cite{schmerl_lat_86} and \cite{schmerl_finite_lat93} and is not easy.
It is an open question whether Theorem \ref{varlt} holds for all \rs\ models.

In another direction, it is observed in \cite[Lemma 7.2]{ks12} that for all $\mc{N}$ and $\mc{K}$, if $\mc{K}\cof\mc{N}$, then  $\Lt_0(\mc{N}/\mc{K})$ is interpretable in $(\mc{N},K)$. This follows from the fact that cofinal Skolem closures over $\mc{K}$ of single elements of $N$ have a particularly simple definition. For all $a$, $b$ in $N$, $a\in \mc{K}(b)$ if and only if $a=(u)_b$  for some $u\in \mc{K}$.\footnote{This notation refers to arithmetic coding of finite sequences: $(u)_x$ is the $x$-th term of the sequence coded by $u$.} Let $m\in K$ be such that $b<m$ and suppose that for a Skolem term $t(x)$, $a=t(b)$. Let $u\in K$ be such that for all $x<m$, $(u)_x=t(x)$. Then $(u)_b=a$. This shows that the relation  $\{\anglebracket{x, y}:\mc{K}(x)\preccurlyeq \mc{K}(y)\}$ is definable in $(\mc{N},K)$.  It follows that if for $i=1,2$, $\mc{K}_i\cof\mc{N}_i$, and $(\mc{N}_1,{K}_1)$ and $(\mc{N}_2,K_2)$ are elementarily equivalent,  then $\Lt_0(\mc{N}_1/{K}_1)$ and $\Lt_0(\mc{N}_2/K_2)$ are elementarily equivalent, and it was asked in \cite{ks12} if in this statement $\Lt_0$ can be replaced by $\Lt$. It turns out that it can not. A counterexample, in which $\Lt_0(\mc{N}_1,\mc{K}_1)$ is an $(\omega+1)$-chain and $(\mc{N}_2,\mc{K}_2)$ is a \rs\ pair elementarily equivalent to $(\mc{N}_1, {K}_1)$, is given in \cite{schmerl_almost_min}.

It follows from Theorem \ref{varlt} and the main result of \cite{ks_london} that if $\mc{M}$ and $\mc{N}$ are \ct\ \as\ models of $\PA$, then 
\[\Lt(\mc{M})\cong\Lt(\mc{N}) \textup{ if and only if }\Aut(\mc{M})\cong\Aut(\mc{N}).\eqno{(*)}\]
In the proof of Theorem \ref{varlt}, it is shown that the left part of the equivalence in $(*)$ is equivalent to $\SSy(\mc{M})=\SSy({N})$, and in the proof of the main theorem of \cite{ks_london} the corresponding equivalence for the automorphism groups is shown to hold for the right part of $(*)$. This leads to the following problem.
\begin{prob}
In the case of \ct\ \as\   models of $\PA$ (or perhaps in greater generality), can any of the two directions of the equivalence $\Lt(\mc{M})\cong\Lt(\mc{N})$ iff $\Aut(\mc{M})\cong\Aut(\mc{N})$ be proven directly, without any reference to the standard systems of the models?
\end{prob}

\section{Open Questions}

There are many questions which remain open about the lattice problem. The following selection is due to Schmerl via private communication. All relevant definitions can be found in \cite[Chapter 4]{ks}.

\begin{enumerate}
    \item Is every finite lattice $L$ which can appear as an interstructure lattice the congruence lattice of some finite algebra? That is, if $L$ is a finite lattice for which there are $\mc{M} \prec \mc{N}$ such that $\Lt(\mc{N} / \mc{M}) \cong L$, is $L$ a congruence lattice of a finite algebra?
    \item Is there a finite lattice $L$ for which there are countable, nonstandard $\mc{M}_0$ and $\mc{M}_1$ for which there is $\mc{N}_0 \succ \mc{M}_0$ such that $\Lt(\mc{N}_0 / \mc{M}_0) \cong L$, but for no $\mc{N}_1 \succ \mc{M}_1$ is $\Lt(\mc{N}_1 / \mc{M}_1) \cong L$?
    \item It is known that if $L$ is a finite lattice that is the congruence lattice of a finite algebra, then so is its dual $L^d$. Does the same result hold for finite lattices $L$ such that there are $\mc{M} \prec \mc{N}$ where $\Lt(\mc{N} / \mc{M}) \cong L$? That is, if $L$ is a finite interstructure lattice, is its dual also an interstructure lattice?
    \item Every finite lattice in the variety generated by $\m3$ is the congruence lattice of a finite algebra. For every  $\aleph_0$-algebraic lattice $L$ in this variety and every countable, nonstandard $\mc{M}$, is there $\mc{N} \succ \mc{M}$ such that $\Lt(\mc{N} / \mc{M}) \cong L$?
    \item For countable $\mc{M}$, what are the possible $\mathfrak{X} \subseteq \mc{P}(\mc{M})$ for which there is an elementary end extension $\mc{N}$ such that $\Lt(\mc{N} / \mc{M}) \cong \n5$ and $\Cod(\mc{N} / \mc{M}) = \mathfrak{X}$?
    \item Is the set of finite lattices for which given (some, all, or a specific) countable, nonstandard $\mc{M}$, there is $\mc{N} \succ \mc{M}$ such that $\Lt(\mc{N}/\mc{M}) \cong L$ computable?
\end{enumerate}
\subsection*{Acknowledgement} We want to thank the referee for a thorough report on this paper. All comments and suggestions in the report  were very helpful in preparation of the final version.

\bibliographystyle{plain}
\bibliography{refs}

\end{document}